\documentclass[english]{amsart}

\usepackage{amsmath}
\usepackage{amssymb}
\usepackage{amsthm}
\usepackage{amscd}
\usepackage{babel}
\usepackage[latin1]{inputenc}
\usepackage{graphicx}

\renewcommand{\subjclassname}{AMS \textup{2000} Mathematics Subject
Classification:\ }

\newtheorem{teor}{Theorem}
\newtheorem{cor}{Corollary}
\newtheorem{prop}{Proposition}
\newtheorem{lem}{Lemma}
\theoremstyle{definition}
\newtheorem{defi}{Definition}

\theoremstyle{definition}
\newtheorem{rem}{Remark}

\newtheorem*{alg}{Algorithm}

\author{Jos\'{e} Mar\'{i}a Grau}
\address{Departamento de Matemáticas, Universidad de Oviedo\\ Avda. Calvo Sotelo, s/n, 33007 Oviedo, Spain}
\email{grau@uniovi.es}
\author{Antonio M. Oller-Marc\'{e}n}
\address{Departamento de Matem\'{a}ticas, Universidad de Zaragoza\\
C/Pedro Cerbuna 12, 50009 Zaragoza (Espa\~{n}a)} 
\email{oller@unizar.es}
\title{A primality test for $Kp^n+1$ numbers}
\begin{document}
\maketitle

\begin{abstract}
In this paper we generalize the classical Proth's theorem for integers of the form $N=Kp^n+1$. For these families, we present a primality test whose computational complexity is $\widetilde{O}(\log^2(N))$ and, what is more important, that requires only one modular exponentiation similar to that of Fermat's test. Consequently, the presented test improves the most often used one, derived from Pocklington's theorem, which usually requires the computation of several modular exponentiations together with some GCD's.
 
\end{abstract}
\subjclassname{11Y11,11Y16,11A51,11B99}

\section{Introduction}
In 1877 P. Pepin (see \cite{PEP}) presented the following result about the primality of Fermat numbers:
%The underlying idea of Proth's theorem was previously presented by P?pin (see \cite{PEP}) in 1877 with the following result about the primality of Fermat numbers.
\begin{teor}[Pepin, 1877]
Let $F_n$ be the $n$-th Fermat number; i.e., $F_n=2^{2^{n}+1}$ with $n>1$. Then, $F_n$ is prime if and only if $3^{\frac{F_n-1}{2}}\equiv -1$ (mod $F_n$).
\end{teor}

Although this theorem has not certified the primality of any new Fermat prime (by 1877 the 5 Fermat primes were already known), it is the first result which leads to a deterministic primality test requiring only one modular exponentiation similar to that of Fermat's test modulo $N$, thus of $\widetilde{O}(\log^{2} N)$ complexity.

One year after, using the same underlying ideas, Proth proved the following primality criterion for number of the form $N=K2^{n}+1$, where $K$ is odd and $K<2^n$ (Proth numbers)

\begin{teor}[Proth, 1878]
Let $N=K2^{n}+1$, where $K$ is odd and $K<2^n$. If $a^{\frac{{N-1}}{2}}\equiv -1$ (mod $N$) for some $a\in\mathbb{Z}$, then $N$ is prime.
\end{teor}

The next important step is the following 1914 result by Pocklington (see \cite{POC}), which is the first generalization of Proth's theorem suitable for numbers of the form $N=Kp^{n}+1$:
 
\begin{teor}[Pocklington, 1914]
Let $N=Kp^{n}+1$ con $K<p^{n}$. If, for some $a\in\mathbb{Z}$:
\begin{itemize}
\item[i)]$a^{N-1}\equiv 1$ (mod $N$)
\item[ii)]$GCD(a^{\frac{N-1}{p}}-1,N)=1$
\end{itemize}
Then, $N$ is prime.
\end{teor}

Proth and Pocklington results are still useful. In fact they are the base of the popular software created by Yves Gallot's (Proth.exe) for the search of Proth and generalized Proth ($N=Kp^{n}+1$) primes.
Other software based in a variation of Pocklington's Theorem presented by Brillhart, Lehmer and Selfridge (see \cite{BRI} or \cite{CRA}) is OpenPFGW with which some records have been broken in different families of integers. For instance, David Broadhurst has recently broken the record for the family $N=2\cdot3^n+1$ (sequence A003306 in the OEIS) certifying primality for $n=1175232$, a number with 560729 digits and the $87$-th biggest known prime (see for instance http://primes.utm.edu/primes/lists/all.txt). An drawback of this software is that it usually requires the use of several bases and, consequently, the computation of several exponentiations modulo $N$.

In recent times the most active researcher looking for primality criteria for numbers of the form $N=Kp^n+1$ has been P. Berrizbeitia. Berrizbeitia and his collaborators have found very efficient criteria for this kind of numbers for a variety of primes $p$ (see \cite{BEBE2,BEBE4,BEIS}). Even though similar criteria had been previously presented by H.C. Williams and his collaborators (see \cite{WIL1,WIL2}), the methodology used by Berrizbeitia et al. shows more clear and efficient. For these generalizations an analogous of Legendre symbol, the \emph{the $m$-th power residue symbol}, has been used. It assumes values over the $m$-th roots of unity and it satisfies a \emph{higher order law of reciprocity}. However, the use of the $m$-th power residue symbol present technical difficulties, mainly due to the fact that the ring $\mathbb{Z}[e^{2\pi i/m}]$ is not a UFD in general. Other authors, such A. Guthmann (see \cite{GUT}) and W. Bosma (see \cite{BOS}), have also given generalizations of Proth's theorem using similar techniques but limited to the case $p=3$.

Our main contribution is a primality criterion for integers of the form $N=Kp^n+1$ with $p$ being any prime and $K<p^n$, using techniques similar to those in \cite{GO} for generalized Cullen Numbers ($N=np^n+1$). These techniques do not require the use of any $m$-th power residue symbol or higher order law of reciprocity. In this way we have achieved an even more clear and efficient methodology than that of Berrizbeitia. In fact, our primality criterion requires only one modular exponentiation $a^{N-1}$ without a previous search of a suitable $a$.

\section{A Generalization of Proth's theorem}

The primality test which follows from Proth's theorem is very useful since, if $N=K2^n+1$ is a prime (Proth Prime), then half the values of $a$ satisfy the condition of the theorem. In particular it is satisfied by those $a$ which are a quadratic non-residue modulo $N$; i.e., such that the Jacobi symbol $(\frac{a}{N})=-1$. This observation is captured in the following version of Proth's theorem:

\begin{teor}[Proth, 1878]
Let $N=K2^{n}+1$, where $K$ is odd and $K<2^n$. Assume that $a\in \mathbb{Z}$ is such that $\left(\frac{a}{N}\right)=-1$, then:
$$N\textrm{is a prime if and only if}\ a^{\frac{{N-1}}{2}}\equiv-1 \textrm{(mod $N$)}.$$
\end{teor}

In spite of the various generalizations presented in the introduction, the most natural generalization of this theorem had not been yet exhibited. We do so in the following result. In what follows $\Phi_p(X)$ will denote the $p$-th cyclotomic polynomial.

\begin{teor}
Let $N=Kp^{n}+1$, where $p$ is a prime, $K<p^n$ and $\gcd(K,p)=1$. Assume that $a\in\mathbb{Z}$ is a $p$-th power non-residue, then:
$$N\textrm{is a prime if and only if}\ \Phi_{p}(\text{a}^{\frac{^{N-1}}{p}})\equiv 0\ \textrm{(mod $N$)}.$$
\end{teor}
\begin{proof}
If $N$ is a prime, then $a^{N-1}\equiv 1$ (mod $N$). Now, $0\equiv a^{N-1}-1=(a^{\frac{N-1}{p}}-1)\Phi_p(a^{\frac{N-1}{p}})$ (mod $N$). Since $a$ is a $p$-th power non-residue, then $a^{\frac{N-1}{p}}-1\not\equiv 0$ (mod $N$) and this implies, $N$ being prime, that $\Phi_p(a^{\frac{N-1}{p}})\equiv 0$ (mod $N$).

Conversely, assume that $\Phi_p(a^{Kp^{n-1}})\equiv 0$ (mod $N$). Put $X=a^K$, then $\Phi_p(X^{p^{n-1}})\equiv 0$ (mod $N$). It follows that $X^{p^n}\equiv 1$ (mod $N$). Now, let $q\leq\sqrt{N}$ be a prime divisor of $N$, then it also holds that $\Phi_p(X^{p^{n-1}})\equiv 0$ (mod $q$) and $X^{p^n}\equiv 1$ (mod $q$). Thus, the order of $X$ in $\mathbb{Z}^{*}_q$ is a divisor of $p^n$, but if $X^{p^j}\equiv 1$ (mod $q$) with $j<n$ would imply that $p=\Phi_p(1)\equiv 0$ (mod $q$) which is clearly a contradiction. Consequently, the order of $X$ in $\mathbb{Z}^{*}_q$ is $p^n$. It follows that $p^n|q-1$ and $p^n<q\leq\sqrt{N}$ and then $p^{2n}\leq N=Kp^n+1$, so $p^n\leq K$ a contradiction.
\end{proof}

The theorem above can be restated in the following way.

\begin{teor}
Let $N=Kp^{n}+1$, where $p$ is a prime and $\gcd(K,p)=1$. If $p^n>K$, then:
$$\Phi_p(a^{\frac{N-1}{p}})\equiv 0\textrm{(mod $N$) $\Leftrightarrow$ $N$ is prime and $a$ is a $p$-th power non-residue modulo $N$}.$$
\end{teor}
\begin{proof}
It is enough to observe that if $\Phi_p(a^{\frac{N-1}{p}})\equiv 0$ (mod $N$), then $N$ is prime (like in the previous proof) and $a\not\equiv x^p$ (mod $N$) for, if it was the case, then $0\equiv\Phi_p(a^{\frac{N-1}{p}})\equiv \Phi_p(x^{N-1})\equiv\Phi_p(1)=p$ (mod $N$); a contradiction.
\end{proof}

This result, like Proth's theorem, is really useful since if $Kp^n+1$ is prime, only $\frac{1}{p}$ of the possible choices for $a$ is a $p$-th power residue modulo $N$. Nevertheless, the interest of this result is mainly theoretical as a genuine generalization of Proth's theorem. An even more useful generalization, not requiring an adequate choice for $a$, will be presented in forthcoming sections.

\section{A generalization of Miller-Rabin primality test}

The so-called Miller-Rabin probabilistic primality test \cite{RAB} test applies to integers in the form $N=K2^n+1$ ($K$ odd) and is based in Fermat's little theorem and in the fact that, the only solutions of $x^2\equiv 1$ (mod $p$) ($p$ prime) are $x\equiv\pm 1$ (mod $p$). In fact we have the following (see \cite[Theorem 3.5.1.]{CRA}): 

\begin{teor} Let $N=K2^{n}+1$ be prime. If $a>1$, then one of the following holds:
\begin{itemize}
\item[i)] $a^{K}\equiv 1$ (mod $N$).
\item[ii)] There exists $0\leq j<n$ such that $(a^{K2^{j}})\equiv -1$ (mod $N$).
\end{itemize}
\end{teor}

 This probabilistic test, in spite of being more demanding than Fermat's test, presents many pseudoprimes (called strong  pseudoprimes) and is specially unreliable if $n$ is small. Nevertheless, for big values of $n$, as in the case of Proth numbers, the test is very reliable and, as we will see in the next section, it allows to certify the primality of the numbers that pass it. 
 
 We must point out that the generalization of Miller-Rabin test is really simple, even though more than two decades passed by until the first publication in this direction. Berrizbeitia and Berry (see \cite{BEBE1}) generalized the Strong Pseudoprime Test introducing the concept $\omega$-prime to base $a$ and more recent work by Berrizbeitia and Olivieri (see \cite{BEOL}) goes in the same direction. Nevertheless, we think that these works do not present a genuine generalization. In fact,  Miller-Rabin test admits a very natural generalization for integers in the form $N=Kp^{n}+1$ with $p$ prime, $K$ even and $\gcd(K,p)=1$. This generalization (that we shall call the $p$-Miller-Rabin test) is based in the following result:

\begin{teor}
Let $p$ be a prime number and $K$ be and even number with $\gcd(K,p)=1$. If $N=Kp^{n}+1$ is prime, then for every integer $a>1$ such that $\gcd(a,N)=1$ one of the following holds:
\begin{itemize}
\item[i)] $a^{K}\equiv1$ (mod $N$).
\item[ii)] There exists $0\leq j\leq n-1$ such that $\Phi_{p}(a^{Kp^{j}})\equiv 0$ (mod $N$).
\end{itemize}
\end{teor}
\begin{proof}
If $N$ is a prime, then $a^{Kp^n}\equiv 1$ (mod $N$). If $a^K\not\equiv 1$ (mod $N$), let $1\leq r\leq n$ be the smallest integer such that $a^{Kp^r}\equiv 1$ (mod $N$). Then $a^{Kp^{r-1}}\not\equiv 1$ (mod $N$) and the primality of $N$ implies that $\Phi_p(a^{Kp^{r-1}})\equiv 0$ (mod $N$) as in Theorem 6. It is enough to put $j=r-1$ to complete the proof.
\end{proof}
  
\begin{defi}
A \emph{$p$-strong probable prime to base $a$} is a number satisfying conditions i) and ii) of Theorem 9 for some $p$, prime divisor of $N-1$. If it is in fact composite, we will say that it is a \emph{$p$-strong pseudoprime to base $a$}.
\end{defi}

This generalization of Miller-Rabin test allows to choose the most appropriate prime factor of $N-1$ in which to base the test. In the case of generalized Proth numbers $N=Kp^{n}+1$ it seems that the prime $p$ should be the most suitable choice; nevertheless, computational experiments reveal
that the number of $q$-strong pseudoprimes does not depend significantly on the chosen divisor of $N-1$. Moreover, the classic
Miller-Rabin test presents in general less pseudoprimes than the proposed generalization. Nonetheless, this new test can be modified to become a
deterministic primality test for Proth numbers ($K<2^n$) and generalized Proth numbers ($N=Kp^n+1$ with $K<p^n$). This modification is the main contribution of this paper and will be developed in the following section.

Also, since $N-1$ will have in general several prime divisors, it makes sense to combine the new test not only using different bases, but
also using different prime divisors of $N-1$. This idea suggests the following definition.

\begin{defi}
A $p$-strong probable prime (resp. $p$-strong pseudoprime) to base $a$ for every $p$ prime divisor of $N-1$, will be denoted as a \emph{complete strong probable prime} (resp. \emph{complete strong pseudoprime}) \emph{to base $a$}.
\end{defi}

Unfortunately, although the concept of complete strong probable prime is more subtle than that of $p$-strong probable prime, computational evidence suggest that it is more convenient to use the test combining different bases rather than different prime divisors of $N-1$. To illustrate this statement it is enough to point out that the smallest 2-strong pseudoprime to bases 2 and 3 is 1373653, while there are 10 complete strong pseudoprimes to base 2 smaller than that number; namely: 2047, 3277, 4033, 8321, 65281, 80581, 85489, 88357, 104653 and 130561.

\section{A sufficient condition for the primality of generalized Proth numbers.}

We will now see that passing the $p$-Miller-Rabin test, together with a bounding condition on $j$ (see Theorem 8), gives a sufficient condition for primality.

\begin{teor}
Let $N=Kp^{n}+1$ where $p$ is a prime and $\gcd(K,p)=1$. If there exists $1\leq j\leq n$ such that:
\begin{itemize}
\item [i)] $\Phi_{p}(2^{Kp^{j-1}})\equiv 0$ (mod $N$).
\item[ii)] $2j> \log_{p}(K)+ n$.
\end{itemize}
Then $N$ is prime.
\end{teor}
\begin{proof}
Put $X=2^{K}$, then $X^{p^{j}}\equiv 1$ (mod $N$). Let $q\leq\sqrt{N}$ be a prime divisor of $N$. It follows that the order of $X$ in $Z_{q}^{*}$ is exactly $p^{j}$. Consequently $p^{j}|q-1$ and $p^{j}<q\leq\sqrt{N}$ from which it follows that $p^{2j}<N=Kp^{n}+1$. Finally, if $p^{2j}\leq Kp^{n}$ then $2j\leq \log_{p}K+n$; a contradiction and the proof is complete.
\end{proof}

\begin{rem}
The theorem above is still true if we replace 2 by any other base $a$. It is enough to put $X=a^K$ in the proof.
\end{rem}

\begin{cor}
Let $N=Kp^{n}+1$ where $p$ is a prime number with $\gcd(K,p)=1$. Let us consider the sequence $S_{0}=2^{K}$, $S_{i}=S_{i-1}^{p}$ for all $i\geq 1$. If for some $j>\frac{1}{2}(\log_{p}(K)+n)$ it holds that $\Phi_{p}(S_{j})\equiv 0$ (mod $N$), then $N$ is prime.
\end{cor}

If we consider the case $p=2$; i.e., the classical Proth numbers, then we get the following corollary.

\begin{cor}
Let $N=K2^{n}+1$ with $K$ an odd integer. Let us consider the sequence $S_{0}=2^{K}$, $S_{i}=S_{i-1}^{2}$ for all $i\geq 1$.
If for some $j>\frac{1}{2}(\log_{2}(K)+ n)$ it holds that $S_{j}\equiv -1$ (mod $N$), then $N$ is prime.
\end{cor}

\section{Algorithm and Computational complexity}

Since 2004, when the polynomial time AKS algorithm was presented (see \cite{AGR}), primality algorithms of general nature were ostracized. That was the case of the deterministic primality test running in $(\log n)^{O(\log\log\log n)}$ time presented by Adleman, Pomerance and Rumely (see \cite{ADL}). This algorithm, later improved by Cohen and Lenstra (see \cite{COH}), is known as the APRCL algorithm. Nevertheless, and despite being one of the cornerstones of Computational Number Theory, AKS algorithm has not been very useful in practice. This is because numbers for which AKS algorithm is faster than the usual ones are beyond current computation capacity. Even the so-called practical versions of the AKS algorithm (see \cite{BER}, for instance) are not fast enough. As a consequence, prime ``hunters'' focus in families of integers for which primality can be determined by useful algorithms. For restricted families of integers much faster algorithms are known, the most celebrated being the Lucas-Lehmer algorithm (see \cite{LUK}), used for Mersenne Numbers, which runs in $\tilde{O}((\log n)^2)$ time. Proth, in \cite{PRO}, gives an algorithm running also in $\tilde{O}(\log n)^2)$ time, which applies to numbers such that $\nu_2(n-1)>\frac12\log_2 n$ where $2^{\nu_2(m)}$ is the biggest power of 2 dividing $m$ and provided an integer $a$ is given such that the Jacobi symbol $\left( \frac{a}{n} \right)=-1$. Proth's algorithm is not deterministic for every $n$. Later, Williams \cite{WIL} or Konyagin and Pomerance \cite{KON} have extended these techniques to wider families of integers.

Unless a surprising discovery is made, the computational complexity of any primality test has a lower bound given by the complexity of the modular exponentiation required by Fermat's test. With this idea in mind, the best that a primality test for an integer $N$ can do is to run in $O(\log^{2}(N) \log(\log(N))\log(\log(\log(N))))$ time. However, even for this complexity, there can be great differences between two different tests depending on the number of modular exponentiations $a^{N-1}$ required. Below we describe an algorithm implementing Corollary 1 which, in fact, requires just one modular exponentiation of the kind $a^{N-1}$ through $n$ modular exponentiations each of them of complexity $O(\log(N) \log(\log(N)) \log(\log (\log(N))))$.

\begin{alg}
\

INPUT: $K,p,n,a$.; $N:=Kp^{n}+1$. S$_{0}:=a^{K}$.

STEP 1: If S$_{0}\equiv1$ (mod $N$)

\ \ \ \ \ \ \ \ \ \ \ \  then RETURN: ``$N$ is a $p$-strong-probable prime to base $a$''. STOP.

STEP 2: For $i=1$ to $n$

\ \ \ \ \ \ \ \ \ \ \ \ \ \ $S_{i}\equiv S_{i-1}^{p}$ (mod $N$)

\ \ \ \ \ \ \ \ \ \ \ \ \ \ If \ \ $S_{i}\equiv 1$ (mod $N$) and $\Phi_{p}
(S_{i-1})\equiv 0$ (mod $N$)

\ \ \ \ \ \ \ \ \ \ \ \ \ \ then Let j:=i. GOTO STEP 3

\ \ \ \ \ \ \ \ \ \ \ \ \ \ If \ \ $S_{i}\equiv 1$ (mod $N$) and $\Phi_{p}
(S_{i-1})\not\equiv 0$ (mod $N$)

\ \ \ \ \ \ \ \ \ \ \ \ \ \ then RETURN: ``$N$ is COMPOSITE'' . STOP

\ \ \ \ \ \ \ \ \ \ \ \ \ \ End

\ \ \ \ \ \ \ \ \ \ \ \ \ \ RETURN: ``$N$ is COMPOSITE''. STOP

STEP 3: If $2j\leq\log_{p}K+n$

\ \ \ \ \ \ \ \ \ \ \ \ \ \ then RETURN: ``$N$ is a $p$-strong-probable prime to base $a$''. STOP.

\ \ \ \ \ \ \ \ \ \ \ \ \ If $2j>\log_{p}K+n$ RETURN: ``$N$ is PRIME''. STOP.
\end{alg}

\begin{prop}
For $N=Kp^n+1$ with fixed $K$ and $p$, the complexity of the algorithm above is $\tilde{O}(\log^{2}(N).$
\end{prop}
\begin{proof}
Only steps 1 and 2 cause complexity, since step 3 is obviously irrelevant.

Complexity of steps 1 is that of the modular exponentiation $a^{K}$ (mod $N$). Taking into account that products modulo $N$ can be performed by Schoenhage-Strassen algorithm (see \cite{SCH}) with complexity:
$$O(\log(N) \log(\log(N))\log(\log(\log(N)))),$$
this is the complexity of step 1.

In step 2 $n$ modular exponentiation with the same complexity as in step 1 are carried out. Thus, since $n=\log_{p}(\frac{N-1}{K})$, the complexity of this step is: $$O(\log^{2}(N) \log(\log(N))\log(\log(\log(N)))).$$

And, summarizing, the whole complexity is $\tilde{O}(\log^{2}(N)$.
\end{proof}

For generalized Proth numbers ($K<p^n$). If we consider $S_{J}:=a^{Kp^J}$ where $J:=\left\lfloor\frac{\log_{p}K+n}{2}\right\rfloor$, it is easy to see that if $S_J\not\equiv 1$ (mod $N$) then the algorithm always certifies the primality or compositeness of $Kp^n+1$. In this case we can consider the following algorithm:

\begin{alg}
\

INPUT: $K,p,n,a$.; $N:=Kp^{n}+1$. $J:=\left\lfloor\frac{\log_{p}K+n}{2}\right\rfloor$. $S_{J}:=a^{Kp^J}$.

\bigskip

STEP 1: If $S_{J}\equiv 1$ (mod $N$)

\ \ \ \ \ \ \ \ \ \ \ \  then RETURN: ``$N$ is a $p$-strong-probable prime to base $a$''. STOP.

\ \ \ \ \ \ \ \ \ \ \ \  else RETURN: ``$N$ will be certified either as prime or composite''.

STEP 2: For $i=J+1$ to $n$

\ \ \ \ \ \ \ \ \ \ \ \ \ \ $S_{i}\equiv S_{i-1}^{p}$ (mod $N$)

\ \ \ \ \ \ \ \ \ \ \ \ \ \ If \ \ $S_{i}\equiv 1$ (mod $N$) and $\Phi_{p}(S_{i-1})\equiv 0$ (mod $N$)

\ \ \ \ \ \ \ \ \ \ \ \ \ \ Then RETURN: ``$N$ is PRIME''

\ \ \ \ \ \ \ \ \ \ \ \ \ \ Else RETURN: ``$N$ is COMPOSITE'' . STOP

\ \ \ \ \ \ \ \ \ \ \ \ \ \ End

\ \ \ \ \ \ \ \ \ \ \ \ \ \ RETURN: ``$N$ is COMPOSITE'' . STOP
\end{alg}

We will now see that for moderately big values of $n$, the probability that the algorithm does not certify the primality of a prime of the form $N=Kp^n+1$ without choosing more that one base is extremely small and that it decreases with $p$. This is not the case for the test based in Pocklington's theorem since, regardless the value of $n$, the use of several bases to certify the primality of $N$ is quite frequent. To do so, we first present a quite well-known lemma.

\begin{lem}
If $N=Kp^n+1$ is prime, the the number of $p^s$-th powers modulo $N$ (different from 0 and 1) is:
$$\frac{N-1}{p^s}-1=Kp^{n-s}-1.$$
\end{lem}

With the use of this lemma we can prove the following proposition.

\begin{prop}
Given a prime $N=Kp^n+1$ ($K<p^n$) and a random base $0<a<n$, the probability that the algorithm returns ``$p$-strong probable prime'' is:
$$\frac{Kp^{\left\lfloor\frac{\log_p(K)+n}{2}\right\rfloor}-1}{Kp^n-1}.$$
\end{prop}
\begin{proof}
The algorithm returns ``$N$ is $p$-strong probable prime'' when $J:=\left\lfloor\frac{\log_p(K)+n}{2}\right\rfloor$ satisfies that $a^J\equiv 1$ (mod $N$). This will happen if $a$ is residual power of order $n-J$ modulo $N$. But, by the previous lemma, the probability that this happens is:
$$\frac{Kp^J-1}{N-2}=\frac{Kp^{\left\lfloor\frac{\log_p(K)+n}{2}\right\rfloor}-1}{Kp^n-1}.$$
\end{proof}

\begin{rem}
For big values of $n$ the probability that a prime of the form $N=Kp^n+1$ is certified as $p$-strong probable prime is about $p^{-n/2}$.
\end{rem}

Steps 1 and 2 in the algorithm perform the computation of the power $a^{N-1}$ (mod $N$) in a controlled way in the sense that if some power $a^{Kp^{i}}\equiv 1$ (mod $N$) the computation stops. Thus, we can say that the computational cost of the algorithm is that of one modular exponentiation of the kind $a^{N-1}$ carried out by $n$ modular exponentiations taking into account that:
$$a^{kp^{n}} =((a^{k})^{p})^{p \stackrel{(n}\cdots p}.$$

Moreover, for values of $p$ with ``many'' 1's or ``many'' 0's in its binary expansion (like for Mersenne or Fermat primes), the presented algorithm can use this fact to perform the $p$-th power in a faster way that with the standard repeat squaring technique; achieving an execution in half the time than the standard modular exponentiation.

To sum up, the presented algorithm improves every primality test requiring more than the computation of a power of the kind $a^{N-1}$ (mod $N$) or similar. It also equals those requiring one such power, even performing better for some particular values of $p$.

\section{Appeal to implementers}
Although the authors have not implemented the proposed algorithm with an appropriate technology, and using Mathematica\textregistered\ only primes up to 100000 digits have been tested, they are in condition to make some considerations that might encourage implementers to create a software based in this paper. Taking into account that our algorithm requires a number of computations similar to that of Fermat's test (or even less) we have compared the time required to certify the primality of the four biggest known primes in the family $N=2\cdot 3^n+1$, recently found by David Broadhurst with the estimated time required by our algorithm. Of course, the runtime of OpenPFGW depends on the ``lucky'' choice of the bases used to perform Pocklington's test (namely, the chosen base $a$ should satisfy $\gcd(a^{\frac{N-1}{P}}-1,N)=1$). OpenPFGW also fails when the tested number is a Fermat pseudoprime for several bases (with a resounding failure when it is a Carmichael number), since it is unable to quickly detect the compositeness of these numbers. However, our algorithm would require only one modular exponentiation of the kind $a^{N-1}$, thus becoming preferable to any other algorithm for generalized Proth numbers. To be true, also our algorithm could require a second choice for the base. But this would happen, for $n=1175232$ with probability about $8.25\times 10^{-280365}$.

In the table below we show the bases used by OpenPFGW (Version 3.4.3) to certify the primality of each $N$ (in one case in needed 7), the runtime in an Intel core2 Duo P7450 @ 2.13 GHz with 4Gb of RAM and the estimated runtime for our algorithm. We also show the ranking of the considered primes among the known primes up to date. All of them are among the 1000 bigger known primes, and the biggest one is among the 100 bigger ones and, remarkably, are among the very few big primes not belonging to the most investigated families: Mersenne, generalized Fermat, Cullen, Woodall, Proth, generalized Cullen and generalized Woodall. It seems to us that the families $Kp^n+1$ have not been deeply investigated except for the case $p=2$.

$$
\begin{tabular}{|l|l|l|l|l|l|l|l|l|l|l|l|l|l|}
\hline
 $N=2\cdot 3^n+1$, $n$=\ &529680&1074726&1086112&1175232\\ \hline  %podemos hacerlo as? $N =$
Number of digits   & 252722 & 512775 &  518208 &560729 \\ \hline
Absolute Ranking   & 895-th & 102-th &  101-th &87-th \\ \hline
%Bases utilizadas for OpenPFGW  & 895-th& 102-th & 101-th & 87-th &  \\ \hline
Bases used by OpenPFGW    & 2,3 & 2,3,17,23,29,31,41& 2  &2,3,5  \\
\hline
Runtime OpenPFGW (in s.)   & 1531. & 21865. & 3220. & 14537 \\ \hline
Estimated runtime our algorithm   & 766. & 3124. & 3220.  & 4845 \\ \hline
\end{tabular}$$

We want to stress the importance of take advantage of the structure of Mersenne and Fermat primes in order to reduce the required time for the modular exponentiations in our algorithm. Consider for instance the search for primes of the form $K\cdot 127^n+1$. Our algorithm requires to perform $n$ modular exponentiations of the kind $b^{127}$. For each of them, performed by the standard repeated squaring algorithm 12 modular products are required, but considering that $b^{127}=b^{128}/b$ only 7 products and a division would be required; a 33\% save. More generally, for $p=2^s-1$ (a Mersenne prime) only $s$ products and a division will be required, while the standard method requires $2(s-1)$ products. Thus, asymptotically, one gets a 50\% save. Moreover, even though $p$ is not a Mersenne or Fermat prime, if there are many 1's or 0's in the binary expansion of $p$ \emph{ad hoc} strategies can be developed in order to optimize the algorithm. This would be the case of primes of the form $2^s\pm 2^t\pm 1$, for instance.

\section*{Acknowledgments}
We are grateful to L. M. Pardo Vasallo for his help with computational complexity aspects. We are also grateful to David Broadhurst, who has helped us to better understand the working of OpenPFGW and whose search for primes of the form $2\cdot3^n+1$ has allowed us to value our algorithm is a more appropriate way.

\bibliography{./refproth}
\bibliographystyle{plain}

\end{document}